\documentclass[12pt,leqno]{article}
\usepackage{amssymb}
\usepackage{amsmath,amsfonts,amssymb}
\numberwithin{equation}{section} \textheight 240mm \textwidth
160mm \topmargin -5mm \oddsidemargin 0mm
\begin{document}

${}$ \vskip 2 cm \centerline{\Large \bf Global solution for the coagulation
equation }

\medskip

\centerline{\Large \bf of water droplets in atmosphere between}

\medskip

\centerline{\Large \bf two horizontal planes.}

\bigskip
\bigskip
\medskip

\centerline{\bf Hanane Belhireche${}^{1)}$, Steave C. Selvaduray${}^{2)}$}
\smallskip

\medskip

\medskip
\centerline{${}^{1)}$ Laboratory of Applied Mathematics and
Modeling,} \centerline{University 8 Mai 1945, Guelma, Algeria,}
\centerline{hanane.belhireche@gmail.com}

\medskip
\centerline{${}^{2)}$ Dipartimento di Matematica, Universit\`a di Torino, Italy,}
\centerline{steave\_selva@yahoo.it, steaveclient.selvaduray@unito.it }

\smallskip
\smallskip
\bigskip
\bigskip

{\small {\bf Abstract.} In this paper we give a global existence
and uniqueness theorem for an initial and boundary value problem
(IBVP) relative to the coagulation equation of water droplets and
we show the convergence of the global solution to the stationary
solution. The coagulation equation is an integro-differential
equation that describes the variation of the density $\sigma$ of
water droplets in the atmosphere. Furthermore, IBVP is considered
on a strip limited by two horizontal planes and its boundary
condition is such that rain fall from the strip. To obtain this
result of global existence of the solution $\sigma$ in the space
of bounded continuous functions, through the method of
characteristics, we assume bounded continuous and small data,
whereas the vector field, besides being bounded continuous, has
$W^{1,\infty}-$ regularity in space.

}
\bigskip
\medskip

\bigskip
\medskip


\bigskip

\vskip 0.8 cm

\section{Introduction}

In \cite{[SF]} has been introduced a model of motion of the air
and the phase transition of water in the three states in the
atmosphere. Since then, many papers have been produced in relation
to this model (see for example \cite{[MF]}, \cite{[BAF2]}
\cite{[AS]}, \cite{[SS1]}, \cite{[SS2]}, etc.). In particular, in
the paper \cite{[MF]}, under some suitable conditions, has been
proved the existence and the uniqueness of the stationary solution
for the equation of water droplets, without taking into account
the condensation or evaporation process, assuming a domain with
two spatial dimensions and a horizontal constant wind. Moreover in
\cite{[BAF2]} has been shown, under weak assumptions, the global
existence and uniqueness of the solution for the same equation,
without taking into account the wind, on a domain with one spatial
dimension. On the other hand, in the paper \cite{[SS2]} has been
proved the existence and the uniqueness, in $W^{1, \infty}$, of
the local solution of an IBVP for the hyperbolic part of the model
introduced in \cite{[SF]},   on a strip bounded by two horizontal
planes with the boundary condition such that rain fall from the
strip.

Now, a significant question which arises, is that to establish the
global existence of the solution of IVBP studied in \cite{[SS2]}.
At the moment, this question is very difficult to treat, therefore, in this paper,
we have decided to study  only one equation of the
hyperbolic part of the model seen in \cite{[SF]}. More precisely,
we consider the coagulation equation of water droplets, including
condensation or evaporation process and effects of wind, in
atmosphere between two horizontal planes (or vertical strip),
supposing known the initial density and the density on the upper
horizontal plane. For this problem, assuming small data and  that
the vertical velocities of droplets are negative,  we prove the
global existence and uniqueness of the solution and we establish
the existence of the stationary solution and the convergence of
the global solution to the stationary solution.

Of course, we know that there exists a well-developed mathematical
theory of coagulation (see for example \cite{[DU]}), but the results
obtained in our paper, although they appear elementary, are not present
in the literature about the coagulation equation.

Afterwards, let us say something about the sections of this paper.
In section 2, we introduce the initial and boundary value problem
for the coagulation equation of water droplets and its stationary
equation associated. In section 3, after having introduced
hypotheses about the velocities of water droplets, using the
method of characteristics, we transform our IBVP in two sets of
infinite Cauchy's problems for ordinary differential equations
(ODEs) and we give the definitions of generalized solutions for
them. Hence, in section 4, we give  the two main theorems of this
paper. In section 5, we study the linearized version of the ODEs
about the coagulation equation obtained in section 3 and we
establish useful estimates to treat our IBVP. In section
6, using some fixed-point arguments, we are able to prove the
first main theorem about the global existence and uniqueness of the solution for the coagulation equation. In
section 7, we prove the second main theorem about stationary
solution.

The authors would like to thank Prof. Hisao Fujita Yashima for having
proposed to us the study of this problem.

\section{Position of the problems.}

\hspace{0,4 cm} We consider the integro-differential equation which describes
the variation of the density $ \sigma (t,x,m) $
of water droplets with mass $m$ in the atmosphere
(see \cite{[SF]}, \cite{[BAF]})
\begin{equation}\label{eq-goutl}
\partial _t \sigma \left( {t,x,m} \right) + \nabla _x
\cdot \left( {\sigma \left( {t,x,m} \right)u
\left( {t,x,m} \right)} \right) + \partial_m
\left( {mh{}_{gl}\left( {t,x,m} \right)\sigma
\left( {t,x,m} \right)} \right) =
\end{equation}
$$
=h_{gl}(t,x,m)\sigma(t,x,m)+\frac{m}{2}\int_{0}^{m}
\beta(m-m',m')\sigma(t,x,m')\sigma(t,x,m-m')dm'+
$$
$$
-m\int_{0}^{\infty}\beta(m,m')\sigma(t,x,m)\sigma(t,x,m')dm'
+g_{0}(m)\big[N_1-\widetilde{N}(\sigma)\big]^{+}(t,x)[Q]^{+}(t,x)+
$$
$$
-g_{1}(m) [Q]^{-}(t,x)\sigma(t,x,m) ,
$$
on the domain
\begin{equation}\label{eq-vapeur}
\mathbb{R}_+ \times \Omega \times \mathbb{R}_+  ,
\end{equation}
where $\Omega=\{ \, x\in \mathbb{R}^3 \, | \,
(x_1 , x_2) \in \mathbb{R}^2 ,
\; 0 < x_3 < 1 \, \} ,$ $t\in \mathbb{R}_+$ is time,
$x = (x_1 , x_2, x_3) \in \Omega $ is a point in three-dimensional
space and $m\in \mathbb{R}_+$ is the mass of a droplet, whereas $[r]^+ $
and $[r]^- $ are the positive and the negative part of $r \in \mathbb{R}$.
In \eqref{eq-goutl}, $u$ and $h_{gl}$
are given functions defined on $\mathbb{R}_+ \times \Omega
\times \mathbb{R}_+ $, $Q $ is a given function on
$\mathbb{R}_+ \times \Omega $, $g_{0}$ and $g_{1}$
are given functions on  $ \mathbb{R}_+$ and $\beta$
is a given function on $\mathbb{R}_+ \times \mathbb{R}_+ $,
whereas $N_1$ and $ \widetilde{N} (\sigma) $ are a positive
constant and a linear functional of $\sigma$ (that we define below) respectively.
Now, we study the equation \eqref{eq-goutl} with the following conditions
\begin{equation}\label{cond-init-sigma}
\sigma (0, x, m ) = \tilde{\sigma}_0 (x, m )
\qquad \mbox{for} \ \ x \in \Omega , \ m \in \mathbb{R}_+ ,
\end{equation}
\begin{equation}\label{cond-lim-sigma}
\sigma (t, x_1, x_2, 1, m ) = \tilde{\sigma}_1 (t, x_1, x_2, m )
\qquad \mbox{for} \ \ t \in \mathbb{R}_+ , \ (x_1, x_2) \in
\mathbb{R}^2, \ m \in \mathbb{R}_+ ,
\end{equation}
where $\tilde{\sigma}_0 $ and $\tilde{\sigma}_1 $ are given
functions defined on $\Omega \times \mathbb{R}_+ $ and
$\mathbb{R}^2 \times \mathbb{R}_+ $ respectively. Moreover,
assuming $ \Gamma_- = \big(\left\{ 0 \right\} \times \Omega
\times \mathbb{R_+} \big) \cup \big( \mathbb{R_+} \times
\mathbb{R}^2  \times \left\{ 1 \right\} \times \mathbb{R_+} \big)
$, we can rewrite \eqref{cond-init-sigma} and
\eqref{cond-lim-sigma} as
\begin{equation}\label{cond-lim-sigma-bis}
\sigma |_{\Gamma _ -  }  = \tilde \sigma  = \left\{ {\begin{array}{*{20}c}
   {\tilde \sigma _0 } & {\mbox{on }\left\{ 0 \right\} \times \Omega  \times \mathbb{R_+}} ,  \\
   {\tilde \sigma _1 } & {\quad \quad \quad \mbox{on }
   \mathbb{R_+} \times \mathbb{R}^2  \times \left\{ 1 \right\} \times \mathbb{R_+}} . \\
\end{array}} \right.
\end{equation}

From a physical point of view, the
function $u(t,x,m) $ is the velocity of the droplet located at $x$ with mass
$m$ and has approximately  the following expression
\begin{equation}\label{def-appr-u}
u(t,x,m) = v(t,x) - \frac{g}{\alpha (m)} e_3 , \quad
e_3 = (0,0,1)^T , \quad (t,x,m) \in \mathbb{R}_+
\times \Omega \times \mathbb{R}_+ ,
\end{equation}
where $ v(t,x) $ is the velocity of air, $g$ is the gravitational
acceleration and $ {\alpha (m)} $ is determined by air friction on
droplet with mass $m$. However, in this paper, we do not assume the expression
\eqref{def-appr-u} for $u$. On the other hand, the function
$h_{gl}(t,x,m)$ is the amount of $H_2O$ that turns from gas to
liquid condensing on a droplet with mass $m$ and $Q = Q(t,x)$ is
the difference between the vapor density $\pi(t,x)$ and the
density $ \overline{\pi}_{vs} (T) $ of a saturated vapor relative
to the liquid state at temperature $T$; more precisely we have
\begin{equation}\label{def-it-Q}
Q(t,x) = \pi(t,x) - \overline{\pi}_{vs} (T(t,x)),
\quad (t,x) \in \mathbb{R}_+ \times \Omega;
\end{equation}
therefore a well approximation for $h_{gl}$ is given by the relation
\begin{equation}\label{def-it-hgl}
h_{gl}(t,x,m)  = \eta (m) Q(t,x),  \quad (t,x,m)
\in \mathbb{R}_+ \times \Omega \times \mathbb{R}_+,
\end{equation}
where the coefficient  $ \eta (m) \geq 0 $ is a lipschitz
function and has a compact support in $\mathbb{R}_+$. Moreover
$\beta(m_1,m_2)$ is the encounter probability between a droplet
with mass $m_1$ and another with masse $m_2$, whereas $g_0 (m)
[N_1 - \widetilde{N} (\sigma) ]^+ $ is the coefficient of
appearance for droplets with mass $m$ and $g_1(m)$ is that of
disappearance for droplets with mass $m$; $N_1$ and $
\widetilde{N} (\sigma) $ are the number of aerosol compared to the
unit of air volume and the number of aerosol already present in droplets
respectively. (For details of physical meaning of these
functions, see \cite{[SF]}, \cite{[BAF]}).

To determine the distribution for $\sigma$, we introduce two numbers
$\overline{m}_a$ and $\overline{m}_A$ (with $0 < \overline{m}_a <
\overline{m}_A < \infty$) and we consider that droplets are absent
apart from an interval $[\overline{m}_a,\overline{m}_A]$, then we have
\begin{equation}\label{def-it-sigma}
\sigma(m) = 0 \qquad \textmd{for } m\in [0,\overline{m}_a[ \cup
]\overline{m}_A, \infty[ ,
\end{equation}
where $\overline{m}_a$, $\overline{m}_A$ correspond respectively to
the lower mass and the upper mass of droplets (see \cite{[PB]}).

Now we suppose that $\beta $ is a continuous function
defined on $\mathbb{R}_+ \times \mathbb{R}_+ $ such that
\begin{equation}\label{cond-beta}
\beta (m_1,m_2) \geq 0 ,\quad \beta (m_1,m_2)= \beta (m_2,m_1)
\ \ \forall (m_1,m_2) \in \mathbb{R}_+ \times \mathbb{R}_+ ,
\end{equation}
\begin{equation}\label{def-C-beta}
\max \big[ \sup\limits_{0<m'<m<\infty } {\frac{m}{2} }
\beta(m -m' ,m') \, , \ \sup\limits_{m,m'\in \mathbb{R}_+ }
{m } \beta(m,m') \big]  < \infty ,
\end{equation}
\begin{equation}\label{cond2-beta}
\beta (m_1, m_2) = 0 \qquad \mbox{if} \ \ m_1 + m_2 \geq
\overline{m}_A  .
\end{equation}
For the functions $g_0 $ and $g_1 $ we suppose that
\begin{equation}\label{cond-g0}
g_0 \in \mathcal{C} ([0, \infty [ ) , \qquad g_0 \geq 0 , \qquad
g_0(m ) = 0 \quad \mbox{if} \ \ m \not \in [ \overline{m}_a ,
\overline{m}_A ],
\end{equation}
\begin{equation}\label{cond-g1}
g_1 \in \mathcal{C} ([0, \infty [ \, )  \quad \mbox{or} \quad g_1
\in \mathcal{C} ( \, ]m_{g1}, \infty [ \, ) , \
\int_{m_{g1}}^{\overline{m}_a+1} g_1(m) dm = \infty , \ \ m_{g1}
\in [0, \overline{m}_a],
\end{equation}
$$
g_1 \geq 0 .
$$
For $ \widetilde{N} (\sigma) $ we assume that
\begin{equation}\label{def-N-tilde}
\widetilde{N} (\sigma)(t,x) = \int_0^\infty n (m) \sigma (t, x, m)
dm,  \quad \forall (t,x) \in \mathbb{R}_+\times \Omega,
\end{equation}
$$
n (\cdot ) \in L^1( \mathbb{R}_+ ) \cap L^\infty
( \mathbb{R}_+ ) , \quad n(m) \geq 0 \quad \forall m \in \mathbb{R}_+ ,
$$
whereas  $N_1$ is a strictly positive constant. 

Now, we consider the stationary equation associated to \eqref{eq-goutl}
\begin{equation}\label{eq-goutl-bis-stat}
 \nabla _x
\cdot \left( {\sigma^{\infty} \left( {x,m} \right)u^*
\left( {x,m} \right)} \right) + \partial_m
\left( {mh{}^*_{gl}\left( {x,m} \right)\sigma^{\infty}
\left( {x,m} \right)} \right) =
\end{equation}
$$
=h^*_{gl}(x,m)\sigma^{\infty}(x,m)+\frac{m}{2}\int_{0}^{m}
\beta(m-m',m')\sigma^{\infty}(x,m')\sigma^{\infty}(x,m-m')dm'+
$$
$$
-m\int_{0}^{\infty}\beta(m,m')\sigma^{\infty}(x,m)\sigma^{\infty}(x,m')dm'
+g_{0}(m)\big[N_1-\widetilde{N}(\sigma^{\infty})\big]^{+}(x)[Q^*]^{+}(x)+
$$
$$
-g_{1}(m) [Q^*]^{-}(x)\sigma^{\infty}(x,m) , \quad (x,m) \in \Omega \times \mathbb{R}_+,
$$
with the following conditions
\begin{equation}\label{cond-lim-sigma-st}
\sigma^{\infty} (x_1, x_2, 1, m ) = \tilde{\sigma}^*_1 (x_1, x_2, m )
\qquad \mbox{for} \ \  \ (x_1, x_2) \in
\mathbb{R}^2, \ m \in \mathbb{R}_+ ,
\end{equation}
where  $\tilde{\sigma}^*_1 $ is a given function defined on
$\mathbb{R}^2 \times \mathbb{R}_+ $ and $u^\ast$, $h^\ast_{gl}$
are the functions defined as in \eqref{def-appr-u},
\eqref{def-it-hgl} respectively and independent on $t$.

Afterwards, we suppose that
\begin{equation}\label{contn-Q-st}
 Q^* \in W^{1,\infty }(\Omega  ), \quad h^*_{gl}(x,m)= \eta(m) Q^*(x) \quad \mbox{for } (x,m) \in \Omega \times \mathbb{R}_+.
\end{equation}

Finally, we conclude this section rewriting the equations \eqref{eq-goutl}  and \eqref{eq-goutl-bis-stat}  respectively as
\begin{equation}\label{eq-goutl.2}
\partial_{t}\sigma+\widetilde{U}\cdot\nabla_{(x,m)}\sigma=
-\widetilde{g}\sigma + \Phi [\sigma ] - \sigma f[\sigma] + h[\sigma]
\end{equation}
and
\begin{equation}\label{eq-goutl.2-st}
\widetilde{U}^*\cdot\nabla_{(x,m)}\sigma^{\infty}=
-\widetilde{g}^*\sigma^{\infty} + \Phi [\sigma^{\infty} ] - \sigma^{\infty} f[\sigma^{\infty}] + h^*[\sigma^{\infty}] ,
\end{equation}
where
\begin{equation}\label{def-nabla}
\nabla_{(x,m)} = \big( \partial_{x_1} , \partial_{x_2} ,
\partial_{x_3} , \partial_{m} \big)^T ,
\end{equation}
\begin{equation}\label{def-U-tilde}
\widetilde{U} (t,x,m) = (u_1(t,x,m),
u_2(t,x,m), u_3(t,x,m), m h_{gl}(t,x,m) ),
\end{equation}
\begin{equation}\label{def-g}
\widetilde{g}(t,x,m) = \nabla_x \cdot u(t,x,m) + \partial_m
(mh_{gl}(t,x,m)) - h_{gl}(t,x,m)+g_{1}(m) [Q (t,x) ]^{-}  ,
\end{equation}
\begin{equation}\label{def-Phi-sigma}
\Phi [ \sigma ] (t,x,m) = \frac{m}{2}
\int_{0}^{m}\beta(m-m',m') \sigma(t,x,m') \sigma(t,x,m-m')dm' ,
\end{equation}
\begin{equation}\label{def-f-sigma}
f [ \sigma ] (t,x,m) =m\int_{0}^{\infty}\beta(m,m')
\sigma(t,x,m')dm' ,
\end{equation}
\begin{equation}\label{def-h}
h [ \sigma ] (t,x,m) = g_{0}(m)\big[N_1-
\widetilde{N}(\sigma)(t,x) \big]^{+} [Q (t,x) ]^{+} ,
\end{equation}
\begin{equation}\label{def-U-tilde-st}
\widetilde{U}^* (x,m) = (u^*_1(x,m),
u^*_2(x,m), u^*_3(x,m), m h^*_{gl}(x,m) ),
\end{equation}
\begin{equation}\label{def-g-st}
\widetilde{g}^*(x,m) = \nabla_x \cdot u^*(x,m) + \partial_m
(mh^*_{gl}(x,m)) - h^*_{gl}(x,m)+ g_{1}(m) [Q^*(x) ]^{-} ,
\end{equation}
\begin{equation}\label{def-h-st}
h^* [ \sigma^{\infty} ] (x,m) = g_{0}(m)\big[N_1-
\widetilde{N}(\sigma^{\infty})(x) \big]^{+} [Q^* (x) ]^{+} ,
\end{equation}
with $(t,x,m) \in \mathbb{R}_+ \times \Omega \times \mathbb{R}_+$.

\section{Assumptions and generalized solutions.}

We make the following assumptions on the velocity $u$  and $Q$
\begin{equation}\label{contn-u}
u \in \mathcal{C}_b ( \mathbb{R}_+ \times \Omega \times
\mathbb{R}_+ ; \mathbb{R}^3) \cap L^1_{loc} \left(
{\mathbb{R}_+;W^{1,\infty } (\Omega  \times
\mathbb{R}_+;\mathbb{R}^3)} \right)  ,
\end{equation}
$$
u_3 \in L_{x_3 }^1 \big(
{0,1;W_{(t_{loc},x_1 ,x_2 ,m)}^{1,\infty } (\mathbb{R}_+ \times
\mathbb{R}^2  \times \mathbb{R}_+)} \big), \quad \nabla \cdot u \in \mathcal{C}_b ( \mathbb{R}_+ \times \Omega
\times \mathbb{R}_+ ),
$$
where $\mathcal{C}_b(X;\mathbb{R})$ is the space of continuous and
bounded functions on a generic metric space $X$,
\begin{equation}\label{def-functional-space-1}
 L^1_{loc} \left(
{\mathbb{R}_+;W^{1,\infty } (\Omega  \times
\mathbb{R}_+;\mathbb{R}^3)} \right) = \bigcap\nolimits_{T > 0} {L^1 \big((0,T);W^{1,\infty } (\Omega  \times
\mathbb{R}_+;\mathbb{R}^3)\big)}
  ,
\end{equation}
and
\begin{equation}\label{def-functional-space-2}
L_{x_3 }^1 \big(
{0,1;W_{(t_{loc},x_1 ,x_2 ,m)}^{1,\infty } (\mathbb{R}_+ \times
\mathbb{R}^2  \times \mathbb{R}_+)} \big)=\bigcap\nolimits_{T > 0} {L^1 \big((0,1);W^{1,\infty } ((0,T) \times
\mathbb{R}^2  \times \mathbb{R}_+)\big)}
;
\end{equation}
moreover there exists a strictly positive constant $A_0$ such that
\begin{equation}\label{borne-u3}
u_3 (t,x,m) \leq - A_0,  \qquad \forall (t,x,m) \in \mathbb{R}_+
\times \Omega \times \mathbb{R}_+
\end{equation}
and we suppose that
\begin{equation}\label{contn-Q}
Q \in  L^1_{loc} \left( {\mathbb{R}_+; W^{1,\infty }
(\Omega  )} \right).
\end{equation}
Therefore the condition \eqref{borne-u3} determines the fall of rain from the strip.

Afterwards, for data $\tilde{\sigma}_0$ and $\tilde{\sigma}_1$, we suppose that
\begin{equation}\label{cont-sigma-bar-0}
\tilde{\sigma}_0  \in \mathcal{C}_b ( \Omega \times \mathbb{R}_+ )
, \qquad \tilde{\sigma}_0 \geq 0 ,
\end{equation}
\begin{equation}\label{cont-sigma-bar-1}
\tilde{\sigma}_1  \in \mathcal{C}_b(\mathbb{R}_+ \times
\mathbb{R}^2 \times \mathbb{R}_+) , \qquad \tilde{\sigma}_1 \geq 0
,
\end{equation}
\begin{equation}\label{cond-sigma-m-bar-aA}
\tilde{\sigma}_0 (x, m ) = 0 , \quad
\tilde{\sigma}_1 (t, x_1, x_2, m ) = 0 \quad \mbox{if} \ \
(t,x) \in \mathbb{R}_+ \times \Omega,\quad m
\not \in [ \overline{m}_a , \overline{m}_A ].
\end{equation}

Furthermore, we assume that
\begin{equation}\label{dis-key-1}
IK \big(1-e^{-J} \big)<J ,
\end{equation}
\begin{equation}\label{dis-key-2}
 \left\| {\tilde \sigma } \right\|_{L^\infty
 \left( {\Gamma _ -  } \right)} <
 \displaystyle\frac{I J \ e^{-J}}{2 \big(J+K I (1-e^{-J}) \big)} ,
\end{equation}
where
\begin{equation}\label{def-I-1}
I = \left\| {\tilde \sigma } \right\|_{L^\infty
\left( {\Gamma _ -  } \right)}  + (1/A_0)\left\| {g_0 }
\right\|_{L^\infty  \left( {\mathbb{R}_+} \right)} N_1
\left\| Q \right\|_{L^\infty  \left( {\mathbb{R}_+
\times \Omega    } \right)},
\end{equation}
\begin{equation}\label{def-J-1}
J = (1/A_0) \big(\| \nabla \cdot u \|_{L^\infty (\mathbb{R}_+
\times \Omega  \times R_ + )} + \|\partial_m (mh_{gl})\|_{L^\infty
(\mathbb{R}_+ \times \Omega \times R_ +)}+
\end{equation}
$$
\ \ \ \ +\| {g_1 }\|_{L^\infty ( \mathbb{R}_+ )}\| Q \|_{L^\infty(
{\mathbb{R}_+ \times \Omega    })} + \left\|
{g_0 } \right\|_{L^\infty \left( {\mathbb{R}_+   } \right)}
\left\| n \right\|_{L^1 \left( {\mathbb{R}_+   } \right)} \left\|
Q \right\|_{L^\infty  \left( {\mathbb{R}_+ \times \Omega
  } \right)} \big),
$$
\begin{equation}\label{def-K-1}
K = (1/A_0) \big(\mathop {\sup }\limits_{m \in \mathbb{R}_+  }
\big| {m\int\limits_0^{ + \infty } {\beta \left( {m,m'}
\right)dm'} } \big| + \mathop {\sup }\limits_{m \in \mathbb{R}_+  }
\big| {\frac{m}{2}\int\limits_0^m {\beta \left( {m - m',m'}
\right)dm'} } \big| \big).
\end{equation}
Now, for the stationary speed $u^*$ we make the following assumptions
\begin{equation}\label{contn-u-st}
u^* \in  W^{1,\infty } (\Omega  \times
\mathbb{R}_+;\mathbb{R}^3)  , \quad  \nabla \cdot u^* \in \mathcal{C}_b (  \Omega
\times \mathbb{R}_+ ),
\end{equation}
\begin{equation}\label{borne-u3-st}
u^*_3 (x,m) \leq - A_0  \qquad \forall (x,m) \in \Omega \times
\mathbb{R}_+ .
\end{equation}
Hence, for  $\tilde{\sigma}^*_1$ we assume that
\begin{equation}\label{cont-sigma-bar-1-st}
\tilde{\sigma}^*_1  \in \mathcal{C}_b(
\mathbb{R}^2 \times \mathbb{R}_+) , \qquad \tilde{\sigma}^*_1 \geq 0
,
\end{equation}
\begin{equation}\label{cond-sigma-m-bar-aA-st}
\tilde{\sigma}^*_1 ( x_1, x_2, m ) = 0 \quad \mbox{if} \ \
(x_1,x_2) \in \mathbb{R}^2,\quad m
\not \in [ \overline{m}_a , \overline{m}_A ],
\end{equation}
\begin{equation}\label{dis-key-2-stat}
 \| {\tilde \sigma^*_1 } \|_{L^\infty
(\mathbb{R}^2 \times \mathbb{R}_+)} <
 \displaystyle
 \frac{{I^*} {J^*} \ e^{-{J^*}}}{2 \big({J^*}+K {I^*} (1-e^{-{J^*}}) \big)} ,
\end{equation}
where
\begin{equation}\label{def-I-bis-stat}
{I^*} =\| {\tilde \sigma^*_1 } \|_{L^\infty (\mathbb{R}^2 \times
\mathbb{R}_+)} + (1/A_0)\left\| {g_0 } \right\|_{L^\infty \left(
{\mathbb{R}_+} \right)} N_1 \left\| Q^* \right\|_{L^\infty \left(
{\Omega    } \right)},
\end{equation}
\begin{equation}\label{def-J-bis-stat}
{J^*} = (1/A_0) \big(\| \nabla \cdot u^* \|_{L^\infty (\Omega
\times R_ + )} + \|\partial_m (mh^*_{gl})\|_{L^\infty (\Omega
\times R_ +)}+
\end{equation}
$$
\ \ \ \ +\| {g_1 }\|_{L^\infty ( \mathbb{R}_+ )}\| Q^* \|_{L^\infty(
{\Omega   })} + \left\| {g_0 }
\right\|_{L^\infty \left( {\mathbb{R}_+   } \right)} \left\| n
\right\|_{L^1 \left( {\mathbb{R}_+   } \right)} \left\| Q^*
\right\|_{L^\infty  \left( {\Omega   } \right)}
\big),
$$
\begin{equation}\label{dis-key-1-st}
I^*K \big(1-e^{-J^*} \big)<J^*.
\end{equation}

Now, after introducing key assumptions about the  velocities and data,
we are ready to give the definitions of generalized solutions for the problems
\eqref{eq-goutl} (or \eqref{eq-goutl.2}), \eqref{cond-lim-sigma-bis} and
\eqref{eq-goutl-bis-stat} (or \eqref{eq-goutl.2-st}),  \eqref{cond-lim-sigma-st}.

First of all, we consider the following Cauchy's problem related
to the flow $X$ associated to the vector field $\widetilde{U}$
\begin{equation}\label{eq-diff-tX}
\begin{cases}
\frac{d}{ds} t(s) = 1 , \\
\frac{d}{ds}X(s)=\widetilde{U}(t(s), X(s)) , \\
(t(0),X(0))=(\tilde{t},\tilde{x_1}, \tilde{x_2},
\tilde{x_3}, \tilde{m}) \in \Gamma_- ,
\end{cases}
\end{equation}
where $X(s)= \big(X_1(s), X_2(s), X_3(s), M(s) \big)$. We observe that
the first equation of \eqref{eq-diff-tX} gives
\begin{equation}\label{eq-diff-solution-t}
t(s) = \tilde{t} + s .
\end{equation}
On the other hand, thanks to the assumptions on  $u$ and $h_{gl}$
(see also the hypothesis \eqref{contn-Q} about $Q$), we deduce there
exists one and only one solution $ X ( \tilde{t} , \tilde{x_1},
\tilde{x_2}, \tilde{x_3},
 \tilde{m};\cdot)$ on $\left[ { - \tilde t, + \infty } \right[
$ for Cauchy's problem \eqref{eq-diff-tX}. Therefore, we
formally can transform the problem \eqref{eq-goutl},
\eqref{cond-lim-sigma-bis} in the following form
\begin{equation}\label{eq1-sigma-cara}
\frac{d}{ds}\sigma(\tilde{t} +s, X (\tilde{t} , \tilde{x_1},
\tilde{x_2}, \tilde{x_3}, \tilde{m}; s)) = - \sigma (\tilde{t} +s,
X (\tilde{t} , \tilde{x_1}, \tilde{x_2}, \tilde{x_3}, \tilde{m}; s)) \times
\end{equation}
$$
\times \big[ \widetilde{g} (\tilde{t} +s, X (\tilde{t} , \tilde{x_1},
\tilde{x_2}, \tilde{x_3}, \tilde{m}; s)) + f [ \sigma ] (\tilde{t}
+s, X (\tilde{t} , \tilde{x_1},
\tilde{x_2}, \tilde{x_3},
\tilde{m}; s)) \big] +
$$
$$
+ \Phi [ \sigma ] (\tilde{t} +s, X (\tilde{t} ,
\tilde{x_1}, \tilde{x_2}, \tilde{x_3},
 \tilde{m}; s )) + h[ \sigma ]  (\tilde{t} +s,
 X (\tilde{t} , \tilde{x_1}, \tilde{x_2}, \tilde{x_3},
 \tilde{m}; s ))  ,
$$
\begin{equation}\label{cond1-sigma-cara}
\sigma(\tilde{t} , \tilde{x_1}, \tilde{x_2}, \tilde{x_3} , \tilde{m}) =
\tilde{\sigma} (\tilde{t} , \tilde{x_1}, \tilde{x_2},
\tilde{x_3},\tilde{m}) \quad  \forall (\tilde{t} , \tilde{x_1},
\tilde{x_2}, \tilde{x_3},\tilde{m}) \in \Gamma_-,
\end{equation}
that is equivalent, in Caratheodory's theory, to the integral equation
\begin{equation}\label{eq1-sigma-cara-intg}
\sigma(\tilde{t} +s, X (\tilde{t} , \tilde{x_1},
\tilde{x_2}, \tilde{x_3}, \tilde{m}; s)) = \tilde{\sigma} (\tilde{t} , \tilde{x_1}, \tilde{x_2},
\tilde{x_3},\tilde{m})  -\int\limits_0^s \Big\{
 \sigma (\tilde{t} +r,
X (\tilde{t} , \tilde{x_1}, \tilde{x_2}, \tilde{x_3}, \tilde{m}; r)) \times
\end{equation}
$$
\times \big[ \widetilde{g} (\tilde{t} +r, X (\tilde{t} , \tilde{x_1},
\tilde{x_2}, \tilde{x_3}, \tilde{m}; r)) + f [ \sigma ] (\tilde{t}
+r, X (\tilde{t} , \tilde{x_1},
\tilde{x_2}, \tilde{x_3},
\tilde{m}; r)) \big] +
$$
$$
+ \Phi [ \sigma ] (\tilde{t} +r, X (\tilde{t} ,
\tilde{x_1}, \tilde{x_2}, \tilde{x_3},
 \tilde{m}; r )) + h[ \sigma ]  (\tilde{t} +r,
 X (\tilde{t} , \tilde{x_1}, \tilde{x_2}, \tilde{x_3},
 \tilde{m}; r )) \Big\}dr .
$$
Hence, we give the following definition

\medskip

{\bf  Definition 3.1.} {\em A continuous solution $\sigma$ for the integral
equation \eqref{eq1-sigma-cara-intg} is  called a generalized solution for the
problem \eqref{eq-goutl} and \eqref{cond-lim-sigma-bis}.
}

\medskip

Now, in a similar way, we can transform the stationary problem
\eqref{eq-goutl-bis-stat}-\eqref{cond-lim-sigma-st} in the following form
\begin{equation}\label{eq1-sigma-cara-bis-stat}
\frac{d}{ds}\sigma^{\infty}( X^* ( \tilde{x_1},
\tilde{x_2}, \tilde{x_3}, \tilde{m}; s)) = - \sigma^{\infty} (
X^* (\tilde{x_1}, \tilde{x_2}, \tilde{x_3}, \tilde{m}; s)) \times
\end{equation}
$$
\times \big[ \widetilde{g}^* ( X^* ( \tilde{x_1},
\tilde{x_2}, \tilde{x_3}, \tilde{m}; s)) + f [ \sigma^{\infty} ] ( X^* ( \tilde{x_1},
\tilde{x_2}, \tilde{x_3},
\tilde{m}; s)) \big] +
$$
$$
+ \Phi [ \sigma^{\infty} ] ( X^* (
\tilde{x_1}, \tilde{x_2}, \tilde{x_3},
 \tilde{m}; s )) + h^*[ \sigma^{\infty} ]  (
 X^* ( \tilde{x_1}, \tilde{x_2}, \tilde{x_3},
 \tilde{m}; s ))  ,
$$
\begin{equation}\label{cond1-sigma-cara-bis-stat}
\sigma^{\infty}( \tilde{x_1}, \tilde{x_2}, 1 , \tilde{m}) =
\tilde{\sigma}^*_1 ( \tilde{x_1}, \tilde{x_2},
\tilde{m}) \quad  \forall ( \tilde{x_1},
\tilde{x_2}, \tilde{m}) \in \mathbb{R}^2 \times \mathbb{R}_+,
\end{equation}
where
the flow $X^* (
\tilde{x_1}, \tilde{x_2}, \tilde{x_3},
 \tilde{m}; s )$ is the solution of the following problem
\begin{equation}\label{eq-diff-tX-bis-stat}
\begin{cases}
\frac{d}{ds}X^*(s)=\widetilde{U}^*( X^*(s)) , \\
X(0)=(\tilde{x_1}, \tilde{x_2},
1, \tilde{m})   ,
\end{cases}
\end{equation}
with $X^*(s)= \big(X^*_1(s), X^*_2(s), X^*_3(s), M^*(s) \big)$.

Of course, the problem
\eqref{eq1-sigma-cara-bis-stat}-\eqref{cond1-sigma-cara-bis-stat}
is equivalent, in Caratheodory's theory, to the following
integral equation
\begin{equation}\label{eq1-sigma-cara-intg-bis-stat}
\sigma^{\infty}( X^* ( \tilde{x_1},
\tilde{x_2}, \tilde{x_3}, \tilde{m}; s)) = \tilde{\sigma}^*_1 ( \tilde{x_1}, \tilde{x_2},
\tilde{m})  -\int\limits_0^s \Big\{
 \sigma^{\infty} (
X^* ( \tilde{x_1}, \tilde{x_2}, \tilde{x_3}, \tilde{m}; r)) \times
\end{equation}
$$
\times \big[ \widetilde{g}^* ( X^* (\tilde{x_1},
\tilde{x_2}, \tilde{x_3}, \tilde{m}; r)) + f [ \sigma^{\infty} ] ( X^* ( \tilde{x_1},
\tilde{x_2}, \tilde{x_3},
\tilde{m}; r)) \big] +
$$
$$
+ \Phi [ \sigma^{\infty} ] ( X^* (
\tilde{x_1}, \tilde{x_2}, \tilde{x_3},
 \tilde{m}; r )) + h^*[ \sigma^{\infty} ]  (
 X ^*(\tilde{x_1}, \tilde{x_2}, \tilde{x_3},
 \tilde{m}; r )) \Big\}dr .
$$

\medskip
Afterwards, we give the following definition of a generalized
solution for the stationary equation
\medskip

{\bf  Definition 3.2.} {\em A continuous solution
$\sigma^{\infty}$ for the integral equation
\eqref{eq1-sigma-cara-intg-bis-stat} is called a generalized
solution for the problem \eqref{eq-goutl-bis-stat} and
\eqref{cond-lim-sigma-st}.  }

\medskip

\section{Main results.}

Now, we are ready to give the first important result of this paper.
\medskip

{\bf  Theorem 4.1.} {\em We assume all hypotheses stated above
on all functions involved in the problem \eqref{eq-goutl} and
\eqref{cond-lim-sigma-bis}. Then there exists one and only one
generalized solution $\sigma$ for the problem \eqref{eq-goutl}, \eqref{cond-lim-sigma-bis}, such that
\begin{equation}\label{prop-sigma}
\sigma \in \mathcal{C}_b (\mathbb{R}_+ \times \Omega \times
\mathbb{R}_+), \quad \sigma \geq 0, \quad \mbox{supp } \sigma
\subseteq \mathbb{R}_+ \times \Omega  \times \left[ {0,\overline m
_B } \right],
\end{equation}
where
\begin{equation}\label{def-mb}
\overline m _B  = \overline m _A \exp \big(
{\left( {1/A_0 } \right)\left\| {h_{gl} } \right\|_\infty  } \big).
\end{equation}
Furthermore  the solution $\sigma$ verifies the inequality
\begin{equation}\label{stima-sigma-5-3-1}
\left\| \sigma  \right\|_{L^\infty  \left( {\mathbb{R}_+ \times \Omega
  \times \mathbb{R}_+ } \right)} \le
   \displaystyle\frac{I J \ e^{-J}}{J+K I (1-e^{-J})}  .
\end{equation}
$($for  $I$, $J$, $K$, see \eqref{def-I-1}-\eqref{def-K-1} $)$.

}

\medskip

Finally, we can give the last main theorem of this paper.

\medskip

{\bf  Theorem 4.2.} {\em We assume all hypotheses stated above
on all functions involved in the problem \eqref{eq-goutl-bis-stat} and
\eqref{cond-lim-sigma-st}. Then there exists one and only one
generalized solution for the
problem \eqref{eq-goutl-bis-stat} and \eqref{cond-lim-sigma-st}, such that
\begin{equation}\label{prop-sigma-bis-stat}
\sigma^{\infty} \in \mathcal{C}_b ( \Omega \times
\mathbb{R}_+), \quad \sigma^{\infty} \geq 0, \quad \mbox{supp } \sigma^{\infty}
\subseteq  \Omega  \times \left[ {0,\overline m
_B } \right],
\end{equation}
where
\begin{equation}\label{def-mb-bis-stat}
\overline m _B  = \overline m _A \exp \big( {\left( {1/A_0 }
\right)\left\| {h^*_{gl} } \right\|_\infty  } \big).
\end{equation}
Furthermore  the solution $\sigma^{\infty}$ verifies the inequality
\begin{equation}\label{stima-sigma-5-3-1-bis-stat}
\left\| \sigma^{\infty}  \right\|_{L^\infty  \left( { \Omega
  \times \mathbb{R}_+ } \right)} \le
   \displaystyle\frac{I^* J^* \ e^{-J^*}}{J^*+K I^* (1-e^{-J^*})}  .
\end{equation}
$($for  $I^*$, $J^*$, see
\eqref{def-I-bis-stat}-\eqref{def-J-bis-stat} $)$.

Finally, if we assume
\begin{equation}\label{supporti-stat}
g_1(m ) = 0, \quad \eta(m)=0 \quad \mbox{if} \ \ m \not \in [ \overline{m}_a ,
\overline{m}_B ],
\end{equation}
\begin{equation}\label{conv-in-l-infinito}
\tilde{\sigma}_1 \mathop  \to \limits_{L^\infty  (\Omega \times
\mathbb{R}_+)} {\tilde{\sigma}^\ast}_1, \qquad Q (t, \cdot
)\mathop \to \limits_{L^\infty  (\Omega )} Q^* ( \cdot ),
\end{equation}
\begin{equation}\label{conv-it-u-to-u*}
u (t, \cdot )\mathop  \to \limits_{L^\infty  (\Omega \times \mathbb{R}_+;\mathbb{R}^3 )} u^* ( \cdot ), \quad \nabla_x \cdot u (t, \cdot
)\mathop  \to \limits_{L^\infty  (\Omega \times \mathbb{R}_+ )}
\nabla_x \cdot u^* ( \cdot ), \quad \mbox{for } t\mathop  \to \infty,
\end{equation}
then
\begin{equation}\label{conv-slz-glob-slz-staz}
\sigma (t, \cdot )\mathop  \to \limits_{L^\infty  (\Omega \times
\mathbb{R}_+)} \sigma ^{\infty} ( \cdot ), \quad \mbox{for } t\mathop  \to \infty,
\end{equation}
where $\sigma$ is defined in Theorem 4.1.
}

\medskip
\medskip
\medskip

\section{Linearization.}

\hspace{0,4 cm} First of all, we observe that from the second
equation of Cauchy's problem \eqref{eq-diff-tX} follows
\begin{equation}\label{def-eq-m-z}
\frac{{dm}}{{dx_3 }} = \frac{{mh_{gl} }}{{u_3 }}
\le m\frac{{\left\| {h_{gl} } \right\|_{L^\infty
( {\mathbb{R}_+ \times \Omega  \times \mathbb{R}_+} )}}}{{A_0 }}.
\end{equation}
If, we define $\overline{m}_B$ such that
\begin{equation}\label{def-mB}
\int\limits_{\overline m _A }^{\overline m _B }
{\frac{{dm}}{m} = } \int\limits_0^1 {\frac{{\left\| {h_{gl} }
\right\|_{L^\infty ( {\mathbb{R}_+ \times \Omega
\times \mathbb{R}_+} )} }}{{A_0 }}dx_3 },
\end{equation}
then we deduce  \eqref{def-mb}. After defining
$\overline{m}_B$, we consider the following cone $K_+$ of
$\mathcal{C}_b(\mathbb{R}_+ \times \Omega  \times \mathbb{R}_+)$
\begin{equation}\label{def-K-}
K_+=\left\{ {\lambda \ | \ \lambda  \in \mathcal{C}_b(\mathbb{R}_+
\times \Omega  \times \mathbb{R}_+), \ \lambda  \ge 0,\ \mbox{supp
} \lambda  \subseteq \mathbb{R}_+ \times \Omega  \times \left[
{0,\overline m _B } \right] } \right\}.
\end{equation}
In this section we study the  linear differential equation
\begin{equation}\label{lin1-eqdiff}
\frac{d}{ds}\sigma(\tilde{t} +s, X (\tilde{t} , \tilde{x_1},
\tilde{x_2}, \tilde{x_3}, \tilde{m}; s)) = - \sigma (\tilde{t} +s,
X (\tilde{t} , \tilde{x_1}, \tilde{x_2}, \tilde{x_3}, \tilde{m}; s)) \times
\end{equation}
$$
\times \Big[ \widetilde{g} (\tilde{t} +s, X (\tilde{t} , \tilde{x_1},
\tilde{x_2}, \tilde{x_3}, \tilde{m}; s))  + f [ \overline{\sigma} ]  (\tilde{t} +s,
X (\tilde{t} , \tilde{x_1}, \tilde{x_2}, \tilde{x_3}, \tilde{m}; s)) \Big] +
$$
$$
+ \Phi [ \overline{\sigma} ] (\tilde{t} +s,
X (\tilde{t} , \tilde{x_1}, \tilde{x_2}, \tilde{x_3},
 \tilde{m}; s )) + h[ \overline{\sigma} ]
 (\tilde{t} +s, X (\tilde{t} , \tilde{x_1}, \tilde{x_2}, \tilde{x_3},
 \tilde{m}; s ))  ,
$$
where $\overline{\sigma} \in K_+$. A first result about this ODE is the following

\medskip

{\bf Lemma 5.1.} \ {\it The equation \eqref{lin1-eqdiff}
with the condition \eqref{cond1-sigma-cara}
has one and only one solution
$\sigma(\tilde{t} +\cdot, X (\tilde{t} , \tilde{x_1},
\tilde{x_2}, \tilde{x_3}, \tilde{m}; \cdot))$
on the interval $\left[ {0,\overline s _1 } \right]$, where
\begin{equation}\label{annulla-x3}
X_3 (\tilde{t} , \tilde{x_1},
\tilde{x_2}, \tilde{x_3}, \tilde{m}; \overline s _1 )=0
\end{equation}
and $\overline s _1$ satisfies the inequality
\begin{equation}\label{borne-s-bar-1}
0 < \overline{s}_1 \leq \frac{1}{A_0} .
\end{equation}
Moreover we have
\begin{equation}\label{posi1-sigma-lin}
\sigma(\tilde{t} +s, X (\tilde{t} , \tilde{x_1},
\tilde{x_2}, \tilde{x_3}, \tilde{m}; s)) \geq 0
\qquad \forall s \in [0 , \overline{s}_1 ] ,
\end{equation}
\begin{equation}\label{supp1-sigma-lin}
\sigma(\tilde{t} +s, X (\tilde{t} , \tilde{x_1},
\tilde{x_2}, \tilde{x_3}, \tilde{m}; s)) = 0
\quad \ \mbox{if} \ \ M(\tilde{t} , \tilde{x_1},
\tilde{x_2}, \tilde{x_3}, \tilde{m};
\overline s _1 ) > \overline{m}_B .
\end{equation}
}

\medskip

{\bf Proof.} \ As $X (\tilde{t} , \tilde{x_1},
\tilde{x_2}, \tilde{x_3}, \tilde{m}; \cdot)$ is well defined
and the equation \eqref{lin1-eqdiff} is linear, the existence and
the uniqueness of the solution for the problem
\eqref{lin1-eqdiff}, \eqref{cond1-sigma-cara}
on $\left[ {0,\overline s _1 } \right]$  follow from the classical theory.
To determine \eqref{borne-s-bar-1}, it is sufficient to  observe that
\begin{equation}\label{4.5.1-s1-1}
\frac{d}{{ds}}X_3 (\tilde{t} , \tilde{x_1},
\tilde{x_2}, \tilde{x_3}, \tilde{m}; s) =
u_3 (\tilde{t} +s, X (\tilde{t} , \tilde{x_1},
\tilde{x_2}, \tilde{x_3}, \tilde{m}; s)) \le -A_0.
\end{equation}
Therefore, the solution $\sigma(\tilde{t} + \cdot, X (\tilde{t} , \tilde{x_1},
\tilde{x_2}, \tilde{x_3}, \tilde{m}; \cdot))$
of the problem \eqref{lin1-eqdiff}, \eqref{cond1-sigma-cara}
can be so represented on the interval
$\left[ {0,\overline s _1 } \right]$ by the expression
\begin{equation}\label{4.5.1-sigma-2}
\sigma(\tilde{t} + s, X (\tilde{t} , \tilde{x_1},
\tilde{x_2}, \tilde{x_3}, \tilde{m}; s)) =
\tilde{\sigma}(\tilde{t}, X (\tilde{t} , \tilde{x_1},
\tilde{x_2}, \tilde{x_3}, \tilde{m}) \times
\end{equation}
$$
\times \exp \Big\{ \int_0^s \big( \widetilde{g}+
f [ \overline{\sigma} ]\big) \big(\tilde{t} +
s', X (\tilde{t} , \tilde{x_1},
\tilde{x_2}, \tilde{x_3}, \tilde{m}; s')\big) ds'\Big\} +
$$
$$
+ \int_0^s \Big[\big(
 \Phi [ \overline{\sigma} ] +
h[ \overline{\sigma} ] \big)(\tilde{t} +
s', X (\tilde{t} , \tilde{x_1},
\tilde{x_2}, \tilde{x_3}, \tilde{m}; s')) \times
$$
$$
\times \exp \Big( \int_{s'}^s \widetilde{g}(\tilde{t} + s'', X (\tilde{t} , \tilde{x_1},
\tilde{x_2}, \tilde{x_3}, \tilde{m}; s'')) ds''\Big) \Big]ds' .
$$
Consequently \eqref{posi1-sigma-lin} can be directly
obtained from \eqref{4.5.1-sigma-2}.

Now to prove \eqref{supp1-sigma-lin} it is sufficient to show that

a) $\sigma(\tilde{t} + s, X (\tilde{t} , \tilde{x_1},
\tilde{x_2}, \tilde{x_3}, \tilde{m}; s)) = 0
$ if $\tilde m > \overline{m}_A$;

b) $M (\tilde{t} , \tilde{x_1},
\tilde{x_2}, \tilde{x_3}, \tilde{m}; s) \le \overline{m}_B $
if $\tilde m \le \overline{m}_A$. \\
The relation a) can be deduced from the consideration that
$M (\tilde{t} , \tilde{x_1},
\tilde{x_2}, \tilde{x_3}, \tilde{m}; \cdot)$ is non decreasing
function and that the conditions \eqref{cond-sigma-m-bar-aA},
\eqref{cond2-beta}, \eqref{cond-g0} are satisfied. To obtain b),
we observe that
\begin{equation}\label{4.5.1-x4-3}
\frac{d}{ds} M (\tilde{t} , \tilde{x_1},
\tilde{x_2}, \tilde{x_3}, \tilde{m}; s) \leq
{\left\| {h_{gl} } \right\|_{L^\infty ( {\mathbb{R}_+
\times \Omega  \times \mathbb{R}_+} )}}
M (\tilde{t} , \tilde{x_1},
\tilde{x_2}, \tilde{x_3}, \tilde{m}; s) ,
\end{equation}
therefore we have
\begin{equation}\label{4.5.1-x4-4}
M (\tilde{t} , \tilde{x_1},
\tilde{x_2}, \tilde{x_3}, \tilde{m}; s) \leq
\tilde{m} \exp \big( {s {\left\| {h_{gl} }
\right\|_{L^\infty ( {\mathbb{R}_+ \times \Omega
\times \mathbb{R}_+} )}} }\big) \le \overline{m}_B .
\end{equation}
$\square$

\medskip

Now we are ready to give the following important result

\medskip

{\bf Lemma 5.2.} \ {\it Let us assume
$\overline{\sigma} \in K_+$ and $\sigma$ the solution of the problem
\eqref{lin1-eqdiff}, \eqref{cond1-sigma-cara}.
Then $\sigma$ is a function on $\mathbb{R}_+ \times \Omega  \times
\mathbb{R}_+$ that belongs to $K_+$. }

\medskip

{\bf Proof.} \ As the relations \eqref{posi1-sigma-lin} and
\eqref{supp1-sigma-lin} are already established,
we just have to show that $\sigma$ is a continuous function on
$\mathbb{R}_+ \times \Omega  \times \mathbb{R}_+$. If we denote
by $\tau_-(t,x_1,x_2,x_3,m)$ the minimum time
of existence for  $X(t,x_1,x_2,x_3,m;\cdot)$ (see \cite{[SS2]}),
then we have the following representation formula for $\sigma$
\begin{equation}\label{tempo-minimo-5-4-3}
\sigma \left( {t,x_1 ,x_2 ,x_3 ,m} \right) = \tilde \sigma
\big(\tau_- (t,x_1 ,x_2 ,x_3 ,m),X(t,x_1,x_2 ,x_3 ,m;\tau_-
(t,x_1 ,x_2 ,x_3 ,m)) \big) +
\end{equation}
$$
+\int\limits_{\tau _-   (t,x_1 ,x_2 ,x_3 ,m)}^t
\Big[\sigma\big(-\widetilde{g} - f\left[{\overline \sigma}\right]
\big) + \Phi \left[ {\overline \sigma  } \right] + h\left[
{\overline \sigma  } \right] \Big] \big(s,X(t,x_1 ,x_2 ,x_3
,m;s)\big) ds .
$$
Now, from the hypotheses \eqref{contn-u}-\eqref{contn-Q} follow
that, thanks to Lemma 4.2 in \cite{[SS2]}, $\tau _ -  \in
W_{(t_{loc},x_1,x_2,x_3,m )}^{1,\infty } (\mathbb{R}_+ \times
\Omega  \times \mathbb{R}_+) =\bigcap\nolimits_{T > 0}
{W^{1,\infty } ((0,T) \times \Omega  \times \mathbb{R}_+)} $ and
${\rm X} \circ \tau _ -  \in W_{(t_{loc},x_1,x_2,x_3,m )
)}^{1,\infty } (\mathbb{R}_+ \times \Omega  \times
\mathbb{R}_+)^4$; therefore, remembering that the given functions
appearing in \eqref{tempo-minimo-5-4-3} are continuous, we deduce
that $\sigma$ is continuous. $\square$

\medskip

{\bf Lemma 5.3.} \ {\it Let $\delta$ be a number such that $0 < \delta \leq 1$. We suppose
that there exists a solution $\sigma$ on
$\mathbb{R}_+\times \Omega_\delta \times \mathbb{R}_+$,
where $\Omega_{\delta}=\mathbb{R}^2 \times (1-\delta,1)$,
for the equation \eqref{eq1-sigma-cara}
with the following condition
\begin{equation}\label{cond1-sigma-cara-delta}
\sigma(\tilde{t} , \tilde{x_1}, \tilde{x_2}, \tilde{x_3} , \tilde{m}) =
\tilde{\sigma} (\tilde{t} , \tilde{x_1}, \tilde{x_2},
\tilde{x_3},\tilde{m}) \quad  \forall (\tilde{t} ,
\tilde{x_1}, \tilde{x_2}, \tilde{x_3},\tilde{m}) \in \Gamma_\delta-,
\end{equation}
where $\Gamma_\delta-=\big[ \left\{ 0 \right\}
\times \mathbb{R}^2 \times (1-\delta,1) \times \mathbb{R}_+ \big]
\cup  \big[\mathbb{R}_+ \times \mathbb{R}^2 \times
\left\{ 1 \right\} \times \mathbb{R}_+ \big].
$
Under these hypotheses, we have
\begin{equation}\label{borne1-sigma-5-4-5}
\left\| \sigma  \right\|_{L^\infty  \left( {\mathbb{R}_+
\times \Omega_\delta  \times \mathbb{R}_+ } \right)} \le
\frac{IJ \ e^{-J}}{J+K I(1-e^{-J})}.
\end{equation}
}

\medskip

{\bf Proof.} \ From the equations relative to the flow associated
to the vector field $\tilde U$ (see \eqref{eq-diff-tX}), we have
\begin{equation}\label{change-der-5-4-5}
\frac{d}{d x_3}\sigma\big(s, X (\tilde{t} , \tilde{x}_1,
\tilde{x}_2,\tilde{x}_3, \tilde{m}; s)\big) =\frac{1}{u}_3
\frac{d}{d s}\sigma\big(s, X (\tilde{t} , \tilde{x}_1,
\tilde{x}_2,\tilde{x}_3, \tilde{m}; s)\big) .
\end{equation}
Therefore, performing integration on \eqref{eq1-sigma-cara}
suitably transformed and after some calculations, we obtain
\begin{equation}\label{estimate-first-5-4-5}
\left\| {\sigma(\cdot,x_3,\cdot)  } \right\|_{L^\infty \left( {
\mathbb{R}_+\times \mathbb{R}^2  \times\mathbb{R}_+ } \right)}\leq
{\left\| {\tilde \sigma  } \right\|_{L^\infty \left( {\Gamma _\delta -
} \right)} } +
\end{equation}
$$
+\int\limits_{x_3 }^{1 } {(1/A_0)  \left\|
 {\sigma(\cdot,z,\cdot)  } \right\|_{L^\infty
 \left( { \mathbb{R}_+\times \mathbb{R}^2  \times\mathbb{R}_+ } \right)}
 \big( \left\| \widetilde{g} \right\|_{L^\infty
 \left( { \mathbb{R}_+\times \Omega_\delta  \times\mathbb{R}_+ } \right)}+
 \left\| {f\left[ \sigma  \right](\cdot,z,\cdot)} \right\|_{L^\infty
  \left( { \mathbb{R}_+\times \mathbb{R}^2  \times\mathbb{R}_+ }
   \right)}\big)} dz+
$$
$$
+\int\limits_{x_3 }^{1 } {(1/A_0) \big[\left\|
{\Phi\left[ \sigma  \right](\cdot,z,\cdot)}
\right\|_{L^\infty  \left( { \mathbb{R}_+\times \mathbb{R}^2
\times\mathbb{R}_+ } \right)} +\left\| {h\left[ \sigma  \right]
(\cdot,z,\cdot)} \right\|_{L^\infty
\left( { \mathbb{R}_+\times \mathbb{R}^2  \times\mathbb{R}_+ }
 \right)}  \big] } dz \le
$$
$$
\le I+\int\limits_{x_3 }^{1 } { \big( J\left\|
{ \sigma  (\cdot,z,\cdot)} \right\|_{L^\infty
\left( { \mathbb{R}_+\times \mathbb{R}^2  \times\mathbb{R}_+ }
\right)} +K\left\| { \sigma  (\cdot,z,\cdot)}
\right\|^2_{L^\infty  \left( { \mathbb{R}_+\times \mathbb{R}^2
\times\mathbb{R}_+ } \right)}  \big) } dz .
$$
Finally,
taking into account the condition \eqref{dis-key-1}
and using the techniques of comparison lemma
we deduce \eqref{borne1-sigma-5-4-5}. $\square$

{\bf Remark 5.4.} We observe that the estimate
\eqref{borne1-sigma-5-4-5} does not depending on $\delta$.

\medskip

\section{Proof of the Theorem 4.1.}

\hspace{0,4 cm}

We start to prove that there exists $0<\delta<1$ such that the
equation \eqref{eq1-sigma-cara} on $\Omega _\delta$ with the
condition \eqref{cond1-sigma-cara-delta} has one and only one
solution. For this purpose, we consider the following operator
$\Lambda _\delta :K_\delta \to \mathcal{C}_b \left( {
\mathbb{R}_+\times \Omega _\delta \times\mathbb{R}_+ } \right)$,
where the domain of $K_\delta$ is so defined
\begin{equation}\label{estimate-def-k-delta}
K_\delta   =\Big\{ \overline \sigma   \in \mathcal{C}_b
\left({\mathbb{R}_+ \times \Omega _\delta   \times \mathbb{R}_+}
\right) \ | \ \overline \sigma   \ge 0,\quad \mbox{supp }\overline
\sigma \subseteq \mathbb{R}_+  \times \Omega _\delta   \times
\left[ {0,\overline m _B } \right] ,
\end{equation}
$$
\left\| {\overline \sigma  } \right\|_{L^\infty  \left(
{\mathbb{R}_+\times \Omega_\delta  \times\mathbb{R}_+ }
\right)}\leq \frac{I J \ e^{-J}}{J+K I (1-e^{-J})}  \Big\}
$$
and $\Lambda_\delta$ is the operator that to each
$\overline\sigma \in K_\delta$
associates the solution $\sigma$
to the linear equation \eqref{lin1-eqdiff}
with the condition \eqref{cond1-sigma-cara-delta}.

Now, we determine an upper bound for $\delta$ such that
$\Lambda_\delta(K_\delta) \subseteq K_\delta$. Using for the
equation \eqref{lin1-eqdiff} an analogous method to that seen for
\eqref{eq1-sigma-cara} in lemma 5.3, we deduce
\begin{equation}\label{estimate-sigma-1}
\left\| {\sigma(\cdot,x_3,\cdot)  } \right\|_{L^\infty  \left( {
\mathbb{R}_+\times\mathbb{R}^2  \times\mathbb{R}_+ } \right)}\leq
{\left\| {\tilde \sigma  } \right\|_{L^\infty  \left( {\Gamma _\delta -
} \right)} } + \int\limits_{x_3 }^{1 } {\overline F\big( {\left\|
{\sigma ( \cdot ,z, \cdot )} \right\|_{L^\infty  \left( {
\mathbb{R}_+\times \mathbb{R}^2  \times\mathbb{R}_+ } \right)} }
\big)} dz,
\end{equation}
where $x_3\in \ ]1-\delta, 1[$ and $\overline F(\cdot)^{}$ is
\begin{equation}\label{def-overline-F}
\overline F\big( {\left\| {\sigma ( \cdot ,z, \cdot )}
\right\|_{L^\infty \left( { \mathbb{R}_+\times \mathbb{R}^2
\times\mathbb{R}_+ } \right)} } \big)= (1/A_0) \Big[ \left\|
{\sigma(\cdot,z,\cdot)  } \right\|_{L^\infty\left( {
\mathbb{R}_+\times \mathbb{R}^2  \times\mathbb{R}_+ } \right)}
\times
\end{equation}
$$
\times \big( \left\| \widetilde{g} \right\|_{L^\infty  \left( {
\mathbb{R}_+\times \Omega_\delta  \times\mathbb{R}_+ } \right)}  +
\left\| {f\left[ \overline \sigma
\right](\cdot,z,\cdot)}\right\|_{L^\infty  \left( {
\mathbb{R}_+\times \mathbb{R}^2  \times\mathbb{R}_+ }
\right)}\big) +
$$
$$
+\left\| {\Phi\left[ \overline \sigma  \right]}
\right\|_{L^\infty  \left( { \mathbb{R}_+\times \Omega_\delta
\times\mathbb{R}_+ } \right)} +\left\| {h\left[ \overline \sigma
\right]} \right\|_{L^\infty  \left( { \mathbb{R}_+\times
\Omega_\delta  \times\mathbb{R}_+ } \right)} \Big]
$$
$$
\le C ( 1+ \left\| {\sigma(\cdot,z,\cdot)  } \right\|_{L^\infty
\left( { \mathbb{R}_+\times \mathbb{R}^2  \times\mathbb{R}_+ }
\right)}    ),
$$
where $C$ is a constant not depending on $\delta$ and
$\overline\sigma$. Therefore, using Gronwall's lemma we have
\begin{equation}\label{estimate-sigma-2}
\left\| {\sigma  } \right\|_{L^\infty  \left( { \mathbb{R}_+\times
\Omega_\delta \times\mathbb{R}_+ } \right)}\leq \exp (C \delta)
\big({\left\| {\tilde \sigma  } \right\|_{L^\infty  \left( {\Gamma
_\delta -} \right)} } +C \delta \big).
\end{equation}
Hence, using \eqref{dis-key-2}, we prove that there exists
$0<\delta_1<1$ such that
\begin{equation}\label{estimate-sigma-3}
\left\| {\sigma  } \right\|_{L^\infty  \left( {
\mathbb{R}_+\times\Omega_\delta  \times\mathbb{R}_+ } \right)}\leq
\frac{IJ \ e^{-J}}{J+K I\big(1- e^{-J} \big)} ,
\end{equation}
for every $0<\delta<\delta_1$. Therefore, thanks to Lemma 5.2
and \eqref{estimate-sigma-3} we deduce that
$\Lambda_\delta(K_\delta)  \subseteq K_\delta$ if
$0<\delta<\delta_1$.

We proceed to study $\Lambda_ \delta$ and in particular we want
to establish an upper bound for $\delta$ such that this map is
a contraction. For this purpose, we consider $\overline \sigma_1$,
$\overline \sigma_2 \in K_\delta$ and we define
$\sigma_j= \Lambda_\delta(\overline \sigma_j)$ with $j=1,2$.
Representing the solutions $\sigma_1$, $\sigma_2$ in integral form,
we deduce, after some transformations, the inequality
\begin{equation}\label{estimate-sigma-diff}
\left\| {\sigma_2(\cdot,x_3,\cdot)- \sigma_1(\cdot,x_3,\cdot)}
\right\|_{L^\infty  \left( { \mathbb{R}_+\times \mathbb{R}^2
\times\mathbb{R}_+ } \right)} \leq
\end{equation}
$$
\leq \frac{1}{A_0}\Big[\big(  \left\| \widetilde{g}
\right\|_{L^\infty  \left( { \mathbb{R}_+\times \Omega _\delta
\times\mathbb{R}_+ } \right)}  + \left\| {f\left[ {\overline
\sigma  _2 } \right]} \right\|_{L^\infty  \left( {\mathbb{R}_+
\times \Omega _\delta   \times \mathbb{R}_+ } \right)}
 \big) \times
$$
$$
\times \int\limits_{x_3 }^{1 } {\left\|{\sigma_2(\cdot,z,\cdot)-
\sigma_1(\cdot,z,\cdot) } \right\|_{L^\infty  \left( {
\mathbb{R}_+\times \mathbb{R}^2  \times\mathbb{R}_+ } \right)}}
dz+
$$
$$
+\delta \big (\left\| {\sigma_1} \right\|_{L^\infty  \left(
{\mathbb{R}_+ \times \Omega _\delta   \times \mathbb{R}_+ }
\right)} \left\| {f\left[ {\overline \sigma _2 } \right]-f\left[
{\overline \sigma  _1 } \right]} \right\|_{L^\infty  \left(
{\mathbb{R}_+ \times \Omega _\delta   \times \mathbb{R}_+ }
\right)} +
$$
$$
+\left\| {\Phi\left[ {\overline \sigma  _2 }
\right]-\Phi\left[ {\overline \sigma  _1 } \right]}
\right\|_{L^\infty  \left( {\mathbb{R}_+ \times \Omega _\delta
\times \mathbb{R}_+ } \right)} +\left\| {h\left[ {\overline \sigma _2 } \right]-h\left[
{\overline \sigma  _1 } \right]} \right\|_{L^\infty  \left(
{\mathbb{R}_+ \times \Omega _\delta   \times \mathbb{R}_+ }
\right)}\big) \Big] \le
$$
$$
\le C\delta\left\| {\overline \sigma  _2 -\overline
\sigma  _1 } \right\|_{L^\infty  \left( {\mathbb{R}_+ \times
\Omega _\delta   \times \mathbb{R}_+ } \right)}+
C\int\limits_{x_3 }^{1 } {\left\|{\sigma_2(\cdot,x_3,\cdot)-
\sigma_1(\cdot,z,\cdot) } \right\|_{L^\infty  \left( {
\mathbb{R}_+\times \mathbb{R}^2  \times\mathbb{R}_+ } \right)}}
dz,
$$
where $C$ is a constant not depending on $\delta$,
$\overline\sigma_1$ and $\overline \sigma_2$. Hence, a direct
application of Gronwall's lemma gives
\begin{equation}\label{stima-it-diff-sigma}
\left\| { \sigma  _2 - \sigma _1 } \right\|_{L^\infty  \left(
{\mathbb{R}_+ \times \Omega _\delta   \times \mathbb{R}_+ }
\right)} \le C \delta \exp(C \delta)\left\| {\overline \sigma  _2
-\overline \sigma  _1 } \right\|_{L^\infty  \left( {\mathbb{R}_+
\times \Omega _\delta   \times \mathbb{R}_+ } \right)}.
\end{equation}
Therefore there exists  $0<\delta_2<\delta_1$
such that $\Lambda_\delta$ is a contraction for every
$0<\delta<\delta_2$. Consequently, the map
$\Lambda_\delta$ admits one and only one fixed point if $0<\delta<\delta_2$.

We have thus proved that  the problem \eqref{eq1-sigma-cara},
\eqref{cond1-sigma-cara-delta} has one and only one solution on
$\Omega_\delta$ if we
assume $0<\delta<\delta_2$. Now, taking into account
the estimate \eqref{borne1-sigma-5-4-5}
and using a simple absurd reasoning
about the maximal width of the strip on which the solution is
defined,
we arrive to show that the solution is defined on the
whole strip $\mathbb{R}_+ \times \Omega  \times \mathbb{R}_+$.
Therefore, the main result
has been proved. $\square$

\medskip

\section{Proof of the Theorem 4.2.}

The proof of the existence and the uniqueness of generalized
solution $\sigma^{\infty}$ together with the properties
\eqref{prop-sigma-bis-stat}-\eqref{stima-sigma-5-3-1-bis-stat} can
be showed in analogous way to that seen for the generalized
solution $\sigma$ in Theorem 4.1. Therefore, we only prove the
last part of this theorem, i.e. we show the convergence
$\sigma(t,\cdot) \to \sigma^\infty$, in $L^\infty  (\Omega \times
\mathbb{R}_+ ) $, for $t \to \infty$.

Let $\delta $ be a number such that $\tilde t>\frac{1}{A_0}$. Hence, we have $\sigma(\tilde{t},\tilde{x_1},\tilde{x_2},1,\tilde{m})=
\tilde{\sigma}_1(\tilde{t},\tilde{x_1},\tilde{x_2},\tilde{m})$ in the integral equation
\eqref{eq1-sigma-cara-intg}. Therefore, we can deduce
\begin{equation}\label{eq-caract}
\sigma(\tilde{t}+s,X(\tilde{t},\tilde{x},\tilde{m},s))=
\tilde{\sigma}_1(\tilde{t},\tilde{x_1},\tilde{x_2},\tilde{m})+
\end{equation}
$$
-\int_0^s \Big(
\big(f[\sigma](\tilde{t}+s',X(\tilde{t},\tilde{x},\tilde{m},s'))
+\tilde{g}(\tilde{t}+s',X(\tilde{t},\tilde{x},\tilde{m},s'))\big)
\sigma(\tilde{t}+s',X(\tilde{t},\tilde{x},\tilde{m},s'))+
$$
$$
+
\big(\Phi [\sigma]+h[\sigma]
\big)(\tilde{t}+s',X(\tilde{t},\tilde{x},\tilde{m},s'))\Big)ds', \quad \mbox{if } \tilde t>\frac{1}{A_0}.
$$

On the other hand, assuming $X(0)=X^\ast(0)$ and recalling the
definitions of the characteristics
$X(\tilde{t},\tilde{x},\tilde{m},s)$ and
$X^\ast(\tilde{x},\tilde{m},s)$ given in \eqref{eq-diff-tX} and
\eqref{eq-diff-tX-bis-stat}, we obtain the following estimate
\begin{equation}\label{eq-it-diff-caract}
|X(\tilde{t},\tilde{x},\tilde{m},s)-X^\ast(\tilde{x},\tilde{m},s)|\leq
\end{equation}
$$
\leq \int_0^s \Big(\big|\tilde{U}(\tilde{t}+s',X(\tilde t, \tilde
x, \tilde m,s'))-\tilde{U}^\ast(X(\tilde t,\tilde x, \tilde m,s'))
\big| +\big|\tilde{U}^\ast(X(\tilde t,\tilde x, \tilde
m,s'))-\tilde{U}^\ast(X^\ast(\tilde x, \tilde m,s')) \big| \Big)
ds' .
$$
Remembering \eqref{contn-Q-st}, \eqref{contn-u-st}, we have that there exists a
constant $C>0$ such that
\begin{equation}\label{eq-it-diff-caract-bis}
|X(\tilde{t},\tilde{x},\tilde{m},s)-X^\ast(\tilde{x},\tilde{m},s)|
\leq \int_0^s \Big(\|u(\tilde{t}+s',\cdot)-u^\ast(\cdot)
\|_{L^\infty (\Omega  \times \mathbb{R}_+; \mathbb{R}^3)} +
\end{equation}
$$
+C\|Q(\tilde{t}+s',\cdot)-Q^\ast(\cdot) \|_{L^\infty (\Omega)} +C
\big|X(\tilde t, \tilde x, \tilde m,s')-X^\ast(\tilde x, \tilde
m,s') \big| \Big) ds'.
$$
Now, applying the comparison lemma to \eqref{eq-it-diff-caract-bis}
and using \eqref{conv-it-u-to-u*}, we obtain the following useful convergence
\begin{equation}\label{conv-it-caract}
|X(\tilde{t},\tilde{x},\tilde{m},s)-X^\ast(\tilde{x},\tilde{m},s)|
\rightarrow_{\tilde{t}\rightarrow \infty} 0.
\end{equation}
Afterwards, a direct application of \eqref{conv-it-caract} into \eqref{eq-caract} gives
\begin{equation}\label{eq-caract-2}
\sigma(\tilde{t}+s,X^\ast(\tilde{x},\tilde{m},s))=
\tilde{\sigma}_1(\tilde{t},\tilde{x_1},\tilde{x_2},\tilde{m}) +
\end{equation}
$$
-\int_0^s \Big(
\big(f[\sigma](\tilde{t}+s',X^\ast(\tilde{x},\tilde{m},s'))
+\tilde{g}(\tilde{t}+s',X^\ast(\tilde{x},\tilde{m},s'))\big)
\sigma(\tilde{t}+s',X^\ast(\tilde{x},\tilde{m},s'))+
$$
$$
+ \big(\Phi [\sigma]+h[\sigma]
\big)(\tilde{t}+s',X^\ast(\tilde{x},\tilde{m},s'))\Big)ds' , \quad
\mbox{if } \tilde t>\frac{1}{A_0}.
$$
Hence, making the difference between the equation
\eqref{eq-caract-2} and the equation
\eqref{eq1-sigma-cara-intg-bis-stat}, after simple algebraic
manipulations, we obtain
$$
|\sigma(\tilde{t}+s,X^\ast(\tilde{x},\tilde{m},s))
-\sigma^\infty(X^\ast(\tilde{x},\tilde{m},s))|\leq
\|\tilde\sigma_1(\tilde{t},\cdot)
-\tilde{\sigma}^\ast_1(\cdot)\|_{L^\infty(\Omega\times\mathbb{R_+})}+
$$
$$
+\int_0^s \Big[\Big(\|\nabla \cdot\big(u(\tilde{t}+s',\cdot)-
u^\ast(\cdot)\big)\|_{L^\infty(\Omega \times \mathbb{R_+}
)}+C'\|Q(\tilde{t}+s',\cdot) -Q^\ast(\cdot)\|_{L^\infty(\Omega)}+
$$
$$
+ C' \|\sigma(\tilde{t}+s',\cdot)
-\sigma^\infty(\cdot)\|_{L^\infty(\Omega\times\mathbb{R_+})}
\Big)\|\sigma(\tilde{t}+s',\cdot)\|_{L^\infty(\Omega\times\mathbb{R_+})}+N_1\|Q(\tilde{t}+s',\cdot)
-Q^\ast(\cdot)\|_{L^\infty(\Omega)}+
$$
$$
+\Big(\|\nabla\cdot u^\ast(\cdot)\|_{L^\infty(\Omega)}+
C'\big(\|Q^\ast(\cdot)\|_{L^\infty(\Omega\times\mathbb{R_+})}+\|\sigma^\infty
(\cdot)\|_{L^\infty(\Omega\times\mathbb{R_+})}+\|\sigma
(\tilde{t}+s',\cdot)\|_{L^\infty(\Omega\times\mathbb{R_+})}\big)\Big)
\times
$$
$$
\times \|\sigma(\tilde{t}+s',\cdot)
-\sigma^\infty(\cdot)\|_{L^\infty(\Omega\times\mathbb{R_+})}\Big]ds'
,
$$
where $C'$ is a positive constant. Now, remembering
\eqref{conv-in-l-infinito} and applying the comparison lemma, we
immediately prove the convergence \eqref{conv-slz-glob-slz-staz}.
$\square$

\medskip
\medskip
\medskip
\medskip



\begin{thebibliography}{15}




\smallskip


\bibitem{[AS]} Ascoli D., Selvaduray Steave C.:
\textit{Wellposedness in the Lipschitz class for a hyperbolic
system arising from a model of the atmosphere including water
phase transitions}. Nonlinear Differ. Equ. Appl., vol 21,
pp. 263-287, 2014.

\smallskip


\bibitem{[BAF]} Belhireche H., Aissaoui M.Z.,
Fujita Yashima H.: \textit{Equations monodimensionnelles
du mouvement de l'air avec la transition de
phase de l'eau.}
Sci. Techn. Univ.
Constantine - A, vol. {\bf 31}, pp. 9--17, 2011.


\smallskip


\bibitem{[BAF2]} Belhireche H., Aissaoui M.Z.,
Fujita Yashima H.: \textit{ Solution globale de l'\'equation de
coagulation des gouttelettes en chute}. 
Ren. Sem. Mat. Univ. Polit. Torino, vol.70, pp.
261--278, 2012.

\smallskip


\bibitem{[BS]} Benzoni Gavage S., Serre D.:
\textit{Multi-dimensional Hyperbolic Partial
Differential Equations.}
Oxford science publications, 2007.

\smallskip


\bibitem{[BR1]} Bressan A., Serre D., Williams M., Zumbrun K., Marcati P.:
\textit{Hyperbolic systems of balance laws.}
 LNM 1911, Springer, 2003.

\smallskip


\bibitem{[BFY]} Buccellato S., Fujita Yashima H.:
\textit{Syst\`eme d'\'equations d'un gaz visqueux
mod\'elisant l'atmosph\`ere avec la force de
Coriolis et la stabilit\'e de l'\'etat d'\'equilibre.}
 Ann. Univ. Ferrara - Sez. VII - Sc. Mat,
vol. 49, pp. 127--159, 2003.

\medskip



\bibitem{[CM]} Chorin A. J., Marsden J. E.:
{\it A mathematical introduction to fluid mechanics}.
III edition, Springer,  1997.

\medskip




\bibitem{[FV]} Filippov A. E.:
\textit{Differential equations with discontinuous righthand sides.}
 Kluwer Academic Pubblishers, 1988.

\medskip


\bibitem{[DU]} Dubovski P. B.:
\textit{Mathematical Theory of Coagulation.}
 (see dubovski.narod.ru/papers/book1.pdf ), 1993.



\medskip


\bibitem{[FCA]} Fujita Yashima H., Campana V.,
Aissaoui M. Z.: \textit{Syst\`eme d'equations d'un modele
du mouvement de l'air impliquant la transition
de phase de l'eau dans l'atmosphere.}
 Ann. Math. Africaines {\bf 2}, pp. 66--92, 2011.


\medskip


\bibitem{[F]} Fujita Yashima H.:
\textit{Modelaci\'on matem\'atica del movimiento de la
atm\'osfera con la transici\'on de fase de agua.}
 Revista Investigaci\'on Operacional, {\bf 34},
\linebreak n. 2, pp. 93--104, 2013.



\medskip

\bibitem{[J1]} Jeffrey A.:
\textit{Quasilinear hyperbolic systems and waves.}
 Pitman publishing limited, 1976.


\smallskip



\bibitem{[LL]} Landau L. L., Lifchitz E. M.:
{\it M\'ecanique des fluides (Physique
th\'eorique, tome 6)} (russian traduction).
Mir, Moscou, 1989.


\medskip


\bibitem{[LTW]}  J. L. Lions, R. Temam, S. Wang:
\textit{New formulations of the primitive equations of
atmosphere and applications.}  Nonlinearity,
vol. {\bf 5}, pp. 237--288, 1992.


\medskip


\bibitem{[MF]} Merad M., Belhireche H., Fujita Yashima H.:
\textit{Solution Stationnaire de l'equation de coagulation de
gouttelettes en chute avec le vent horizontal.}
Rend. Sem. Mat. Univ. Padova,
vol. 129, pp. 225--244, 2013.


\medskip






\bibitem{[PB]} Prodi F., Battaglia A.:
\textit{Meteorologia - Parte II, Microfisica.}
Grafica Pucci, Roma, 2004. (see also site:
http://www.meteo.uni-bonn.de/mitarbeiter/battaglia/teaching.html ).




\medskip






\bibitem{[SF]} Selvaduray S. C., Fujita Yashima H.:
\textit{Equazioni del moto dell'aria con la transizione
di fase dell'acqua nei tre stati: gassoso, liquido e solido.}
Memorie della classe di scienze Fisiche, Matematiche
e Naturali, Accademia
delle Scienze di Torino, Serie V, volume 35, pp.37--69,   2011. (see also  http://hdl.handle.net/2318/734).


\medskip


\bibitem{[SS1]} Selvaduray S. C.:
\textit{On a quasilinear hyperbolic system relative to an atmospheric model on the transition of water defined on the whole space.}     Journal of Mathematics and System Science, Volume 4, Number 9, September 2014.



\medskip


\bibitem{[SS2]} Selvaduray S. C.:
\textit{An Initial and  boundary value problem on a strip
for a large class of quasilinear hyperbolic systems arising from an atmospheric model.}
(Arxive:1411.2119v1, 2014).

\smallskip
\medskip


\smallskip


\smallskip


\end{thebibliography}
\end{document}